# THE POWERS OF 9 AND RELATED MATHEMATICAL TABLES FROM BABYLON

*Mathieu Ossendrijver (Humboldt University, Berlin)*

Late-Babylonian (LB) mathematics (450-100 BC), represented by some sixty cuneiform tablets from Babylon and Uruk, is incompletely known compared to its abundantly preserved, well-studied Old-Babylonian (OB) predecessor (1800-1600 BC).[1] With the present paper, sixteen fragments from Babylon, probably belonging to 13 different tablets, are added to this corpus. Two remarkable tablets represent a hitherto unknown class of very large factorization tables that can be adequately described as Babylonian examples of number crunching (Section I). Most other fragments belong to tables with reciprocals (II) and squares (III). Finally, two fragments contain multiplications of one kind or another (IV).

All fragments belong to the Babylon collection of the British Museum, which includes thousands of astronomical tablets (diaries, related observational texts and mathematical astronomy) originating from unscientific excavations conducted in the 19th century. It has long been recognized that these astronomical tablets originate from Babylon, where they were written by generations of astronomers. Their dates range from the Neo-Babylonian to the Parthian era (550 BC - 50 AD), but most are Late-Achaemenid or Seleucid (350-150 BC). The same provenience and distribution of dates can be assumed for the mathematical fragments in the Babylon collection, all of which belong to lots with numerous astronomical texts. From the Achaemenid era onwards the astronomers in Babylon are believed to have been employed by the Esagila temple (Beaulieu 2006). An invocation to Bēl and Bēltiya, titulary gods of the Esagila, preserved on fragment A1, supports the assumption that the mathematical texts from Babylon were written by the same group of scholars.

Size and content of the LB mathematical corpus continue to raise questions about the scope and possible applications of mathematics in the LB era and its degree of autonomy from mathematical astronomy. Near 400 BC the astronomers developed mathematical algorithms for predicting astronomical phenomena, building on the OB mathematical apparatus and utilising centuries of observational data compiled in astronomical diaries. Some significant mathematical innovations are attested only in mathematical astronomy.[2] Within the mathematical corpus, important innovations were thus far identified mainly in the few extant problem texts.[3] Most other tablets of this corpus are traditional tables with reciprocals, squares or fourth powers, differing from the OB tables mainly by a denser distribution and greater length of their principle numbers, which contain up to seven sexagesimal digits. While probably fashioned by astronomers, none of these tables are actually needed in mathematical astronomy, at least not for its routine operations.[4] Furthermore, the squares and fourth powers contain up to twenty-five digits, far beyond any practical application. The factorization tables presented in Section I add significant new evidence for these and other supra-utilitarian aspects of LB mathematics.

## Sexagesimal place-value notation, regular numbers and reciprocals

In Babylonian mathematics, numbers are expressed using sexagesimal place-value notation, which operates analogously to our decimal notation. Numbers are represented as sequences of digits, each assuming a value between 0 and 59. Every digit is associated with a power of 60, which decreases from one to the next digit. From the Achaemenid era onwards, vanishing intermediate digits (0) were indicated by a sign consisting of two small, vertically aligned *Winkelhakens*.[5] The same sign, here also transliterated as 0, could be placed between multiples of 10 and a following 1-9 in order to mark that they represent separate digits.[6] As opposed to our decimal system, which is an absolute notation,

1 For a list of published LB mathematical tablets and a discussion of their archival context cf. Robson (2008), pp. 214-262, 337-344.

2 E.g., operations involving quantities of undetermined magnitude ('variables') and a concept of additive and subtractive numbers (Ossendrijver 2012, pp. 29-31).

3 E.g., regarding the nature of the problems, their formulation and solution methods: cf. Friberg (1990); Friberg, J., Hunger, H., al-Rawi, F. (1997); Høyrup (2002), pp. 387-399; Robson (2008); Ossendrijver (2012), pp. 26-27.

4 This was pointed out by Robson (2008), p. 261. However, during the formative period of mathematical astronomy, tables of reciprocals must have been used for converting divisions into multiplications by a reciprocal (Section II).

5 In astronomical tables occasionally also initial vanishing digits. The same sign functions as a separator (:) between words or numbers, e.g. in Texts E and I.

6 In Babylon this convention is less common than in Uruk, but it is used rather systematically in Texts A and B (in Texts C-M no instance of a digit 10-50 followed by a digit 1-9 is preserved).

the Babylonian sexagesimal notation is relative, in the sense that the power of 60 corresponding to each digit is not indicated and must be inferred from the context. Since all tablets presented here are context-free tables, the absolute value of the numbers is not determined. In the transliterations this feature is preserved by separating all digits by periods (.). If the absolute value of a number is known, it can be represented by placing a semicolon (;) between the digit pertaining to 1 and that pertaining to 1/60, and commas between all other digits. For instance, 6;40 (=6+40/60), 0;6,40 ($=6/60+40/60^2$), 6,40,0 ($=6\cdot 60^2+40\cdot 60$), etc, are absolute readings of 6.40.

An important role is played by so-called regular numbers, which are numbers, say *s*, of which the reciprocal, 1/*s*, exists as a finite sequence of sexagesimal digits. This is true if and only if *s* contains no other factors besides 2, 3 and 5, the prime factors of 60 ($= 2^2\cdot 3\cdot 5$). All regular numbers therefore have the form

$$s = 2^p \cdot 3^q \cdot 5^r,$$

where *p, q* and *r* are integers.[7] Once they are known, 1/*s* is obtained as $(0;30)^p \cdot (0;20)^q \cdot (0;12)^r$, where 0;30=1/2, 0;20=1/3, and 0;12=1/5. In the relative notation 1/*s* is equivalent to $30^p \cdot 20^q \cdot 12^r$, e.g. $1/9 = 20^2 = 6.40$.[8] In all transliterations, reconstructed signs are enclosed by square brackets and shaded gray for additional clarity. Usually an explanatory column is added on the right side. Computer-based tools for generating and querying tables of sexagesimal numbers were essential for identifying and reconstructing all fragments.

## I. Factorization tables

Texts A and B belong to a hitherto unknown class of factorization tables. The numerical virtuosity displayed in these tables is unprecedented; in particular, the numbers in lines Obv. 1'-8' of Text B contain between twenty-six and thirty sexagesimal digits, which makes them the longest numbers found in a cuneiform text and probably in all antiquity. As to the purpose of these tables no definite explanation can be offered yet. It seems unlikely that they were relevant for practical computations in accounting or astronomy, which rarely involve numbers with more than seven digits.[9] Their purpose must be located within the context of scholarly mathematics. As will be argued, they may have served as proof that their initial number had been computed correctly. Alternatively, or simultaneously, they may be viewed as didactive examples that demonstrate number-theoretical regularities relating the final digits of a number to the factor to be used in the factorization algorithm (Friberg 1999). In Text A the final digits alternate between 9 and 21, a diagnostic feature of powers of 9 prompting the use of a factor 6.40. Similarly, the final digits 12, 36, 48 and 24 in Obv. 0'-39' of Text B suggest a factor 5. A third possibility, namely that Texts A and B belong to two-columnar algorithms for computing reciprocals, can be virtually ruled out for Text A and also appears less likely for Text B. By exemplifying unparalleled, supra-utilitarian numerical skills and, perhaps, a notion of numerical proof and an interest in number-theoretical regularities, factorization tables shed significant new light on LB mathematics, for which they imply a greater autonomy from mathematical astronomy then previously thought.

### Text A. Factorization table for $9^{46}$

| | |
|---|---|
| Fragment A1 | BM 34249 (Sp, 356)<br>size: 3.5 x 3.3 x 1.4-1.9 cm<br>content: Obv. 1-7 partly preserved |
| Fragment A2 | BM 32401 (76-11-17, 2134) +32707 (2475)<br>size: 6.6 x 4.7 x 2.3 cm<br>content: Obv. 6-15 partly preserved |
| Fragment A3 | BM 34517 (Sp, 641)<br>size: 3.3 x 4.0 x 1.3 cm<br>copy: *LBAT* 1646<br>content: Obv. 22-28 partly preserved |

7 Note that some combinations of *p, q* and *r* result in the same sequence of digits. This is because a change from (*p, q, r*) to (*p+2n, q+n, r+n*), where *n* is an integer, amounts to a multiplication by $60^n$.

8 In the absolute notation this corresponds to 1/9 = 0;6,40 or 1/9,0 = 0;0,6,40, etc.

9 One unusual astronomical procedure text, Ossendrijver (2012) Nr. 102, contains several numbers with nine digits.

Obverse

| | | |
|---|---|---|
| 1 | *ina a-mat* $^{d}$EN *u* $^{d}$GAŠAN-⸢*$ia_2$*⸣ [*liš-lim*] | |
| 2 | 16.34.3⸢9⸣.[52.40.21.26.52.57.35.56.49.50.37.38.58.13.38.4.44.57.15.3.37.21] | $9^{46}$ |
| 3 | 1.50.31.⸢5⸣.[51.9.2.59.13.3.59.38.52.17.30.59.48.10.53.51.39.41.40.24.9] | $9^{45}$ |
| 4 | 12.16.47.1[9.1.0.19.54.47.6.37.39.8.36.46.38.41.12.39.4.24.37.49.21] | $9^{44}$ |
| 5 | 1.21.51.5[5.26.46.42.12.45.14.4.11.0.57.25.10.57.54.44.20.29.24.12.9] | $9^{43}$ |
| 6 | 9.5.[46.9.38].⸢31.21.25.1.33.47.53.26⸣.[22.47.53.6.4.55.36.36.1.21] | $9^{42}$ |
| 7 | 1.0.[38.27.44.16].⸢4⸣9.2.46.50.25.19.16.15.51.⸢5⸣[9.14.0.32.50.44.0.9] | $9^{41}$ |
| 8 | [6.44.16.24.5]⸢5⸣.12.6.58.32.16.8.48.[28.25.46.34.53.23.38.58.13.21] | $9^{40}$ |
| 9 | [44.55.9.2]⸢6⸣.8.0.46.30.15.7.38.43.9.3[1.50.32.35.57.39.48.9] | $9^{39}$ |
| 10 | [4.59.27.42].⸢5⸣4.13.25.10.01.40.50.58.⸢7⸣.[43.32.16.57.19.44.25.21] | $9^{38}$ |
| 11 | [33.16.24.4]6.1.29.27.46.51.12.⸢1⸣[9.47.31.30.15.13.2.11.36.9] | $9^{37}$ |
| 12 | [3.41.49.2]⸢5⸣.6.49.56.25.⸢12.21⸣.[22.11.56.50.01.41.26.54.37.21] | $9^{36}$ |
| 13 | [24.38.49].⸢27.25.3⸣2.56.⸢8⸣.[2.22.27.59.38.53.31.16.19.24.9] | $9^{35}$ |
| 14 | [2.44.18].49.42.50.19.34.1[3.35.49.46.37.39.16.48.28.49.21] | $9^{34}$ |
| 15 | [18.15].⸢25.31.25.35.30⸣.[28.10.38.51.50.51.1.52.3.12.9] | $9^{33}$ |
| 16 | [2.1.42.50.09.30.36.43.7.50.59.5.39.0.12.27.1.21] | $9^{32}$ |
| 17 | [13.31.25.34.23.24.4.47.32.19.53.57.40.01.23.0.9] | $9^{31}$ |
| 18 | [1.30.09.30.29.16.0.31.56.55.32.39.44.26.49.13.21] | $9^{30}$ |
| 19 | [10.01.3.23.15.6.43.32.59.30.17.44.56.18.48.9] | $9^{29}$ |
| 20 | [1.6.47.2.35.0.44.50.19.56.41.58.19.35.25.21] | $9^{28}$ |
| 21 | [7.25.13.37.13.24.58.55.32.57.59.48.50.36.9] | $9^{27}$ |
| 22 | [49.28.10.48.9.2]⸢6⸣.3⸢2.50⸣.[19.46.38.45.37.21] | $9^{26}$ |
| 23 | [5.29.47.52.1.2].⸢5⸣6.5⸢8⸣.[55.31.50.58.24.9] | $9^{25}$ |
| 24 | [36.38.39.6.46].⸢5⸣9.39.[52.50.12.19.49.21] | $9^{24}$ |
| 25 | [4.4.17.40.45.1]3.17.⸢45⸣.[52.14.42.12.9] | $9^{23}$ |
| 26 | [27.8.37.51.41].28.38.2⸢5⸣.[48.18.1.21] | $9^{22}$ |
| 27 | [3.0.57.32.24.3]⸢6⸣.30.56.⸢1⸣[2.2.0.9] | $9^{21}$ |
| 28 | [20.06.23.36.4].3.2⸢6⸣.[14.40.13.21] | $9^{20}$ |
| 29 | [2.14.2.37.20.27.2.54.57.48.9] | $9^{19}$ |
| 30 | [14.53.37.28.56.20.19.26.25.21] | $9^{18}$ |
| 31 | [1.39.17.29.52.55.35.29.36.9] | $9^{17}$ |
| 32 | [11.1.56.39.12.50.36.37.21] | $9^{16}$ |
| 33 | [1.13.32.57.41.25.37.24.9] | $9^{15}$ |
| 34 | [8.10.19.44.36.10.49.21] | $9^{14}$ |
| 35 | [54.28.51.37.21.12.9] | $9^{13}$ |
| 36 | [6.3.12.24.9.1.21] | $9^{12}$ |
| 37 | [40.21.22.41.0.9] | $9^{11}$ |
| 38 | [4.29.2.31.13.21] | $9^{10}$ |
| 39 | [29.53.36.48.9] | $9^{9}$ |
| 40 | [3.19.17.25.21] | $9^{8}$ |
| 41 | [22.8.36.9] | $9^{7}$ |
| 42 | [2.27.37.21] | $9^{6}$ |
| 43 | [16.24.9] | $9^{5}$ |
| 44 | [1.49.21] | $9^{4}$ |
| 45 | [12.9] | $9^{3}$ |
| 46 | [1.21] | $9^{2}$ |
| 47 | [9] | 9 |
| 48 | [1] | 1 |

Critical commentary:
Obv. 1: 'By the command of Bēl and Bēltiya [may it be succesful]', or: '... [may it (the tablet) remain intact]'.
Obv. 22: The digit 3⸢2⸣ is erroneously copied as 3⸢8⸣ in *LBAT*.

Text A is preserved in three fragments, probably belonging to a single tablet whose dimensions were at least 16 x 18 x 2.3 cm. By a fortunate coincidence A1 preserves the upper left corner, thus providing crucial evidence for the full extent of the table. Its reverse side, which belongs to the lower left corner, is intact but not inscribed. The reverse sides of A2 and A3 are destroyed. As far as known, this is the only LB mathematical tablet from Babylon preserving an invocation (Obv. 1). The mention of Bēl and Bēltiya implies that it was written by scholars associated with the Esagila temple, most likely astronomers. Each of the following lines contains a single number whose final digit is aligned on the right side of the table, resulting in a triangular layout. However, most of the preserved digits are not strictly aligned.[10] On the right side of the transliteration there is a modern representation of each number as a power of 9.

The reconstruction of Text A proceeded from A2, which is mentioned in C.B.F. Walker's unpublished catalogue of astronomical fragments. However, the length of the numbers suggested a mathematical rather than an astronomical content. This prompted a comparison with BM 55557 (Britton 1991-3), a table of regular numbers with up to twenty-five digits, which are fourth powers. After querying a computer-generated table of all regular numbers with up to nine digits and their fourth powers, Obv. 6, 8, 10, 12 and 14 turned out to match $(3^{21})^4 = 9^{42}$, $(3^{20})^4 = 9^{40}$, $(3^{19})^4 = 9^{38}$ and $(3^{18})^4 = 9^{36}$, respectively. It was then obvious that Obv. 7, 9, 11 and 13 should contain the intermediate odd powers of 9, which was confirmed by inspecting a computer-generated table of powers of 9. A subsequent search yielded fragments A2 and A3. Not a single scribal error had to be assumed anywhere in Obv. 2-28. The obvious question, whether the same fragments might not allow a different reconstruction, given the size of the gaps, can be answered with a clear no. No other regular numbers with up to thirty digits, apart from those reconstructed here, match the digits preserved in Obv. 2-15.[11] The resulting uniformly decreasing sequence of powers of 9 leaves little doubt that the reconstruction is correct up to Obv. 28. Below that it is also very probable, being based on the reasonable assumption that the factorization continues down to 1 without scribal errors. It is assumed here that the entire table was written on the obverse but, as suggested by Text B, some part below Obv. 15 may have been written on the reverse.

The closest parallels of Texts A and B are two LB oval-shaped tablets from Uruk with factorization algorithms for computing reciprocals of regular numbers.[12] On these tablets, regular numbers with up to six digits are factorized down to 1 as in Text A. The reciprocal of each number was written in an adjacent column which was filled from the bottom to the top. For example, if we set out from 12.9 = $9^3$ then the algorithm proceeds as shown in Table 1. The restriction to regular numbers applies to both columns: if the initial number is not regular it cannot be factorized down to 1.

| factor | | *s* (col. i) | 1/*s* (col. ii) | | factor |
|---|---|---|---|---|---|
| 6.40 ↓ | $9^3 =$ | 12.9 | 4.56.17.46.40 | $= (6.40)^3$ | |
| 6.40 ↓ | $9^2 =$ | 1.21 | 44.26.40 | $= (6.40)^2$ | 6.40 ↑ |
| 6.40 ↓ | $9^1 =$ | 9 | 6.40 | $= (6.40)^1$ | 6.40 ↑ |
| | $9^0 =$ | 1 | 1 | $= (6.40)^0$ | 6.40 ↑ |

Table 1: Algorithm for computing reciprocals of powers of 9

10 For instance, the preserved digits in Obv. 8 are shifted by about one digit.

11 In fact, A1 alone contains sufficient information for reconstructing the entire table. On A3 each sequence of preserved digits is by itself too short for a unique identification, but their succession makes the identification with decreasing powers of 9 very probable.

12 *SpTU* 4 176 (Friberg 1999; 2000) and *SpTU* 5 316 (Friberg 2007, p. 453). The former, perhaps also the latter, belonged to the library of the Šangû-Ninurta family, which suggests an approximate date 415±30 BC (Clancier 2009, pp. 397, 404; Robson 2008, pp. 227-240). Friberg (1999) interprets *SpTU* 4 176 as a school tablet, but the oval shape, typical for certain OB school tablets, is unknown for LB school texts (Gesche 2000, pp. 18, 43-53). Moreover, their content is beyond the curriculum of the scribal school (Gesche 2000, pp. 136-140). If they were written by students then this can only have happened in the framework of a more advanced level of education.

This sheds light on Text A, by suggesting that there may have been a second column containing the corresponding power of 6.40 for each number in the first column. However, with every multiplication these reciprocals rapidly grow to unreasonable lengths (Table 1). The number $(6.40)^{46}$ would occupy sixty-eight digits, probably too much even for a Babylonian mathematician and requiring a column as wide as 35 cm. Moreover, a search in the Babylon collection did not yield any fragment preserving a portion of this column. I therefore conclude that Text A most likely did not include a second column with reciprocals.

In order to assess the significance of the number equivalent to $9^{46}$ note that it is the fourth power of 2.1.4.8.3.0.27 (= $3^{23}$). In so-called extended tables of reciprocals this is the first principle number ($s$) beyond the standard table, that is, whose initial digit exceeds 1 (cf. Section II). Hence the sexagesimal number equivalent to $9^{46}$ might appear in an extended table of fourth powers, but no such table has been found yet.[13]

Another LB mathematical fragment from Babylon to be mentioned here, *LBAT* 1644, partly preserves a computation of $9^{46}$ as a sexagesimal number by squaring $9^{23}$ through digit-wise multiplication and subsequent addition of the products (Friberg 2007, pp. 456-9).[14] It is well possible that both tablets were written by the same scribe.[15] This suggests that the purpose of Text A may have been to verify that its initial number (Obv. 2), as computed in *LBAT* 1644, is a correct representation of $9^{46}$, by repeatedly multiplying it with 6.40=1/9, until 1 is reached. It is worth mentioning that *LBAT* 1644 ends with a correct representation of $9^{46}$, even though the preceding computation from which it was obtained contains an error. As suggested by Friberg (2007), *LBAT* 1644 may be an incomplete copy of a tablet with a correct version of the algorithm. Alternatively, the outcome of *LBAT* 1644 may have been obtained not by squaring as one is made to believe, but by multiplying 9 with itself 46 times, which is how Text A can be read from bottom to top. Since multiplying by 9 requires only half as many operations as multiplying by 6.40, the factorization in Text A is probably fictitious, in the sense that the numbers were copied in reverse order from an original on which the powers of 9 increased downwards. Both attestations of the number $9^{46}$ might go back to that original computation.

## Text B. Factorization table for $9^{11}\cdot 12^{n}$ ($n \geq 39$)

Fragment B1 BM 42744 (81-7-1, 508) +45977 (81-7-6, 420) +46008 (452)
size: 11.2 x 7.3 x 2.3 cm
content: Obv. 1'-14' partly preserved

Fragment B2 BM 34958 (Sp2, 479)
size: 7.0 x 2.2 x 2.5 cm
copy: *LBAT* 1642
content: Obv. 18'-21' partly preserved

Obverse

| | | |
|---|---|---|
| 0' | [2.5.10.44.5.44.37.31.5.10.56.27.9.29.31.52.8.21.18.15.0.19.36.13.31.57.59.16.13.26.24] | $9^{11}\cdot 12^{40}$ |
| 1' | [10.25.53.40.28.43.7.35.25.54.42.15.47.27.39.20.41.46.31.15.1.38.1.7].⌜39.4⌝[9.56.21.7.12] | $9^{11}\cdot 12^{39}$ |
| 2' | [52.9.28.22.23.35.37.57.9.33.31.18.57.18.16.43.28.52].⌜36.15.8.10⌝.5.38.1⌜9⌝.[9.41.45.3]6 | $9^{11}\cdot 12^{38}$ |
| 3' | [4.20.47.21.51.57.58.9.45.47.47.36.34.46.31.23].⌜37.24.2⌝3.1.15.40.⌜50⌝.[28.11].⌜3⌝5.48.28.48 | $9^{11}\cdot 12^{37}$ |
| 4' | [21.43.56.49.19.49.50.48.48.58.58.2.53.52.36.58].⌜7.x1.55.6.18.24.12.20⌝.[57].59.2.24 | $9^{11}\cdot 12^{36}$ |
| 5' | [1.48.39.44.6.39.9.14.4.4.54.50.14.29.23.4].50.35.9.35.31.32.1.1.4[4].⌜49.5⌝5.12 | $9^{11}\cdot 12^{35}$ |
| 6' | [9.3.18.40.33.15.46.10.20.24.34.11.12.26.55.24.12.5]5.47.57.37.40.5.8.4⌜4⌝.9.36 | $9^{11}\cdot 12^{34}$ |
| 7' | [45.16.33.22.46.18.50.51.42.2.50.56.2.14.37.1].⌜4.3⌝[8].59.48.8.20.25.⌜43⌝.40.48 | $9^{11}\cdot 12^{33}$ |
| 8' | [3.46.22.46.53.51.34.14.18.30.14.14.40.11.13].⌜5⌝.5.2⌜3.14⌝.59.<0>.⌜41⌝.42.8.3⌜8⌝.24 | $9^{11}\cdot 12^{32}$ |
| 9' | [18.51.53.54.29.17.51.11.32.31.11.13.20.56].⌜5⌝.25.26.56.14.55.3.2⌜8.30.4⌝3.12 | $9^{11}\cdot 12^{31}$ |
| 10' | [1.34.19.29.32.26.29.15.57.42.35.56.6.44.40].27.7.14.41.14.35.17.22.[3]⌜3⌝.36 | $9^{11}\cdot 12^{30}$ |
| 11' | [7.51.37.27.42.12.26.19.48.32.59.40.33.43.22.1]⌜5.3⌝6.13.26.12.56.26.⌜5⌝[2.4]8 | $9^{11}\cdot 12^{29}$ |
| 12' | [39.18.7.18.31.2.11.39.2.44.58.22.48.36.51.1]⌜8⌝.1.7.11.4.42.14.2[4] | $9^{11}\cdot 12^{28}$ |
| 13' | [3.16.30.36.32.35.10.58.15.13.44.51.54.3.4.16].⌜30⌝.05.35.55.23.⌜31⌝.1⌜2⌝ | $9^{11}\cdot 12^{27}$ |

13 The two known LB tables of fourth powers, BM 55557 (Britton 1991-3) and BM 32584 (Friberg 2007, p. 455), deal with principle numbers between 1 and 2.

14 BM 34601 (Sp2,76+759) (in Friberg 2007, pp. 458-9, erroneously referred to as MS 36401 and BM 46301).

15 *LBAT* 1644 might belong to the reverse of the same tablet as Text A, but a physical join could not be established.

| | | |
|---|---|---|
| 14' | [16.22.33.2.42.55.54.51.16.8.44.19.30.15.21.22.30.27].⸢59.36⸣.[57.36] | $9^{11} \cdot 12^{26}$ |
| 15' | [1.21.52.45.13.34.39.34.16.20.43.41.37.31.16.46.52.32.19.58.4.48] | $9^{11} \cdot 12^{25}$ |
| 16' | [6.49.23.46.7.53.17.51.21.43.38.28.7.36.23.54.22.41.39.50.24] | $9^{11} \cdot 12^{24}$ |
| 17' | [34.6.58.50.39.26.29.16.48.38.12.20.38.1.59.31.53.28.19.12] | $9^{11} \cdot 12^{23}$ |
| 18' | [2.50.34.54.13.17.12.26.24.3.11.1.43.10.09.57.39.27.21].⸢36⸣ | $9^{11} \cdot 12^{22}$ |
| 19' | [14.12.54.31.6.26.2.12.0.15.55.8.35.50.49].⸢48.1⸣7.16.48 | $9^{11} \cdot 12^{21}$ |
| 20' | [1.11.4.32.35.32.10.11.0.1.19.35.42.59.1]⸢4⸣.9.1.26.24 | $9^{11} \cdot 12^{20}$ |
| 21' | [5.55.22.42.57.40.50.55.0.6.37.58].⸢34⸣.56.10.45.7.12 | $9^{11} \cdot 12^{19}$ |

Reverse

| | | |
|---|---|---|
| 1 | [29.36.53.34.48.24.14.35.0.33.9.52.54.40.53.45.36] | $9^{11} \cdot 12^{18}$ |
| 2 | [2.28.4.27.54.2.1.12.55.2.45.49.24.33.24.28.48] | $9^{11} \cdot 12^{17}$ |
| 3 | [12.20.22.19.30.10.06.4.35.13.49.7.2.47.2.24] | $9^{11} \cdot 12^{16}$ |
| 4 | [1.1.41.51.37.30.50.30.22.56.9.5.35.13.55.12] | $9^{11} \cdot 12^{15}$ |
| 5 | [5.8.29.18.7.34.12.31.54.40.45.27.56.9.36] | $9^{11} \cdot 12^{14}$ |
| 6 | [25.42.26.30.37.51.2.39.33.23.47.19.40.48] | $9^{11} \cdot 12^{13}$ |
| 7 | [2.8.32.12.33.9.15.13.17.46.58.56.38.24] | $9^{11} \cdot 12^{12}$ |
| 8 | [10.42.41.2.45.46.16.6.28.54.54.43.12] | $9^{11} \cdot 12^{11}$ |
| 9 | [53.33.25.13.48.51.20.32.24.34.33.36] | $9^{11} \cdot 12^{10}$ |
| 10 | [4.27.47.6.9.4.16.42.42.2.52.48] | $9^{11} \cdot 12^{9}$ |
| 11 | [22.18.55.30.45.21.23.33.30.14.24] | $9^{11} \cdot 12^{8}$ |
| 12 | [1.51.34.37.33.46.46.57.47.31.12] | $9^{11} \cdot 12^{7}$ |
| 13 | [9.17.53.7.48.53.54.48.57.36] | $9^{11} \cdot 12^{6}$ |
| 14 | [46.29.25.39.4.29.34.4.48] | $9^{11} \cdot 12^{5}$ |
| 15 | [3.52.27.8.15.22.27.50.24] | $9^{11} \cdot 12^{4}$ |
| 16 | [19.22.15.41.16.52.19.12] | $9^{11} \cdot 12^{3}$ |
| 17 | [1.36.51.18.26.24.21.36] | $9^{11} \cdot 12^{2}$ |
| 18 | [8.4.16.32.12.1.48] | $9^{11} \cdot 12$ |
| 19 | [40.21.22.41.0.9] | $9^{11}$ |

Critical commentary:
Obv. 0': Even though nothing remains of this line, there is sufficient clay above Obv. 1' to suggest that it may have been on the tablet with a certain probability.
Obv. 2': There is a large open space between 10 and 5, but no separation mark (0).
Obv. 4': Between the digits 7 and 1 there are traces of a sign similar to 20; nothing is expected.
Obv. 6': There is a large open space between 40 and 5, but no separation mark (0).
Obv. 8': The expected 0 between 59 and 41 is lacking.
Obv. 19': The digit 17 is erroneously copied as 27 in *LBAT.*

Text B is preserved in two fragments probably belonging to a single tablet whose dimensions were at least 22 x 11 x 2.5 cm. Both are inscribed on one side which is probably the obverse; their other side is destroyed. B1 consists of three newly joined fragments. It does not preserve any edge of the tablet, but not much clay may be missing on its right side and above it. B2 belongs to the lower edge of the tablet. Each line contains a single number whose first digit must have been aligned on the left side of the table, which explains the triangular layout. Up to twelve final digits of each number are preserved, but they are not strictly aligned, i.e. the amount of space used for each digit is not always the same. There are two minor scribal errors (Obv. 4', 8'), neither of which was passed on to the numbers below or above it, implying that the tablet was copied from an original without these errors. On the right side of the transliteration there is a modern representation of each number as a product of powers of 9 and 12.

The reconstruction of Text B proceeded from a photograph of BM 46008 found in the *Nachlass* of A. Aaboe and O. Neugebauer which is kept in the Institute for the Study of the Ancient World (ISAW) in New York. After querying a computer-generated table of all regular numbers with up to thirty digits, the well-preserved digits on that fragment (Obv. 5'-13') matched regular numbers with (*p, q, r*) decreasing in steps of (2, 1, 0) from (70, 57, 0) to (54, 49, 0) i.e. *s* decreasing from $9^{11}\cdot12^{35}$ to $9^{11}\cdot12^{27}$. A search in the Babylon collection then yielded fragments B2 and BM 42744. The latter physically joined the top edge of BM 46008, extending Text B upwards by four powers of 12. As it turned out, BM 45977 had been joined to BM 46008, thus confirming Obv. 2'-11', while B2 confirmed that Text B continues until Obv. 21'. Again, no other regular numbers with up to thirty digits, apart from the ones reconstructed here, match the preserved digits in Obv. 2'-13', 21'. The digits preserved in Obv. 19' match sixteen regular numbers with up to thirty digits, those in Obv. 20' match four such numbers. In each case the expected number is among them. The resulting uniformly decreasing sequence of powers of 12 leaves little doubt about the reconstruction of Obv. 1'-21', in spite of the large gaps. There is some uncertainty as to how the table continues above Obv. 1' and below Rev. 19. That Obv. 1' may be preceded by another line containing the number equivalent to $9^{11}\cdot12^{40}$ is suggested by the amount of clay available above Obv. 1'. The table probably continued on the reverse until $9^{11}$ was reached in Rev. 19, followed by the same sequence of 11 decreasing powers of 9 as in Obv. 38-48 of Text A. Alternatively the factorization proceeded from there in steps of 1/3 = 20, in which case twenty-two lines with decreasing powers of 3, from $3^{21}$ to 1, must be reconstructed below Rev. 19.

As was suggested for Text A, Text B may have contained another column, now lost, which provided the reciprocal of each number in the preserved column, i.e. from 1 in the last line to $(6.40)^{11}$ in Rev. 19, $(6.40)^{11}\cdot5$ in Rev. 18, until $(6.40)^{11}\cdot5^{39}$ in Obv. 1'. Since these reciprocals are not prohibitively long - the sexagesimal representation of $(6.40)^{11}\cdot5^{39}$ occupies 32 digits - there is no compelling practical reason for ruling out the presence of this column. The alignment of the partly preserved numbers on their left side (compare Text A) may suggest that they form the right column of a two-columnar computation of the reciprocal of $(6.40)^{11}\cdot5^{39}$. In that interpretation the partly preserved column contains successive multiplications by 12=1/5 to be read from bottom to top. However, a search in the Babylon collection did not yield any portion of the hypothetical missing column, in which $(6.40)^{11}\cdot5^{39}$ should be factorized down to 1. Furthermore, the purpose of Text B may not be to compute reciprocals but, as was suggested for Text A, to prove that its initial number is correct. After repeatedly multiplying it by 5=1/12 and then eleven times by 6.40=1/9 or twenty-two times by 20=1/3, the outcome should be 1.[16] In that interpretation there was no other column. For these reasons I have not included it in the reconstruction. The significance of the initial number, $9^{11}\cdot12^{n}$ ($n\geq39$), remains wholly unclear.

## II. Tables of regular numbers and their reciprocals

Five fragments, Texts C-G, preserve portions of a more or less standardized two-columnar table with approximately one hundred regular numbers between 1 and 2 and their reciprocals. The total number of such tablets has thereby reached fifteen.[17] For a reconstruction of the complete table and the underlying criteria governing the inclusion of the numbers cf. Britton (1991-3) and Friberg (2007, pp. 461-4). The index numbers (1-100) from Friberg's standard table are provided on the right side of each transliteration. None of the fragments contains a pair that is not in that table, but Text G may lack some pairs. In the LB era, the term for reciprocals known from OB mathematics, IGI.x.GAL$_2$.BI y, 'the reciprocal of x is y', was usually abbreviated to IGI x y.[18] Since IGI would be written near the left edge, it is often broken away. With the exception of Text C and Text H, an OB-style table of reciprocals, all tablets are assumed to employ that formulation.

The function of these tables is not completely clear. Very few numbers in the standard table have more than seven digits, about the maximum number of digits encountered in Babylonian mathematical astronomy. Indeed divisions do not occur in mathematical astronomy, since they were reformulated as multiplications by reciprocals, presumably with the help of the present tables.[19] However, it is not clear why nearly all are limited to principle numbers *s* = 1...2. It would appear to be equally important to have access to reciprocals of numbers *s*>2, but such extended tables are very rare and not one of them covers the range *s*>5.[20]

16 One may take the parallel with Text A one step further and propose the existence of a text similar to *LBAT* 1644 in which the initial number of Text B is computed by digit-wise squaring. This would support the idea that X.1' is not the first line, because $9^{11}\cdot12^{n}$ is not the square of a regular number if *n* is odd. However, a search in the Babylon collection did not yield any fragment of a computation of $9^{11}\cdot12^{40}$ by digit-wise squaring of $3^{11}\cdot12^{20}$.

17 Previously published texts: Vaiman (1961), Texts B-F; Britton (1991-3), Texts A-C; Friberg (2005), Nos. 72-76. The latter five fragments are assumed to belong to two different tablets.

18 Cf. Britton (1991-3), Text A.

19 Ossendrijver (2012), pp. 25-26.

20 Only three examples are known: *TU* 31 (Neugebauer 1935-7, pp. 14-22) and *SpTU* 4 174, both from Uruk, contain reciprocals of numbers 1<*s*<3; Aaboe (1965), Text I, a fragment from Babylon, contains reciprocals of numbers 4<*s*<5.

## Text C = BM 36065 (Sp3, 611)

Obverse

| | i = *s* | ii = 1/*s* | Nr |
|---|---|---|---|
| 1 | [1.0].⸢45 | 5⸣9.⸢15⸣.33.⸢20⸣ | 1 |
| 2 | [1.1.2.6].33.⸢45 | 58⸣.58.56.3⸢8.24⸣ | 2 |
| 3 | [1.1].26.24 | 58.35.37.30 | 3 |
| 4 | [1.1].30.33.45 | 58.31.39.35.18.51.6.40 | 4 |
| 5 | [1.1].⸢4⸣3.42.13.20 | ⸢58⸣.19.12 | 5 |
| 6 | [1.2].12.28.48 | 57.52.13.⸢20⸣ | 6 |
| 7 | [1.2.30] | ⸢57⸣.3⸢6⸣ | 7 |
| 8 | [1.3.12.3]⸢5.33⸣.20 | ⸢56.57.11.15⸣ | 8 |
| 9 | [1.3.1]⸢6.52.30 | 56.53.20⸣ | 9 |
| 10 | [1.4] | ⸢56.15⸣ | 10 |
| 11 | [1.4].⸢48 | 55.33⸣.20 | 11 |
| 12 | ⸢1⸣.[5.6].⸢15 | 55.17.45⸣.36 | 12 |
| 13 | [1.5].⸢36⸣.[3]⸢6⸣ | ⸢54.52⸣.10.51.51.6.40 | 13 |
| 14 | [1.5.50.3]⸢7⸣.2.13.20 | ⸢54.40.30⸣ | 14 |
| 15 | [1.5.55.4].41.15 | ⸢54.3⸣[6.4]⸢8⸣ | 15 |
| 16 | [1.6.21.18.4]⸢3.12 | 54⸣.[15.12].⸢30⸣ | 16 |
| 17 | ⸢1⸣.[6].⸢40⸣ | 54 | 17 |
| 18 | ⸢1⸣.[7.30] | ⸢53.20⸣ | 18 |
| 19 | [1.8.16] | ⸢5⸣[2.44].⸢3⸣.[45] | 19 |
| 20 | [1.8.20.37.30 | 52].⸢40.29.37⸣.[46.40] | 20 |
| 21 | [...] | ⸢52⸣.[...] | ? |
| 22-25 | [...] | ⸢50⸣ [...] | ? |

Reverse

| | |
|---|---|
| 1'-4' | (faint illegible traces) |
| 5' | [xxxxxx] ⸢xx 4$^?$7 x⸣ [xx] |
| 6' | [xxxxxx] ⸢x 18$^?$ xxxx 30⸣ [xx] |
| 7' | [xxxxxx] ⸢xx 8.17$^?$.4$^?$.34$^?$ x⸣ [xx] |
| 8' | [xxxxxxxxx] ⸢x 8 xx⸣ [xx] |
| 9' | [xxxxxxxxxxx] ⸢x 5$^?$⸣ [xx] |
| 10' | [xxxxxxxxxx] ⸢xxxxx⸣ [xx] |
| 11' | [xxxxx] ⸢xxxxxx⸣ |
| 12' | [xxxxx] ⸢9$^?$.20 xxx⸣ |
| 13' | [xxxxx] ⸢xxxxx⸣ |
| 14' | [xxxxx] ⸢x 30$^?$ xxxxx 3 x 37.30$^?$⸣ |
| 15' | [xxxxx] ⸢8$^?$.3$^?$2.37$^?$ xx 1.50 xx⸣ |
| 16' | [xxxxx] ⸢2.1.4.8⸣.3.0.2⸢7⸣ |
| 17' | [xxxxx] ⸢10$^?$ x 20$^?$ xx 10$^?$ xx⸣ [x] |
| 18'-20' | (illegible) |

Critical commentary:
Obv. 12, 17, 18: Only faint traces remain of the leftmost sign, probably a 1, less likely IGI, 'reciprocal of'.
Rev. 5'-17': All readings with question marks are very uncertain.
Rev. 15': The 8$^?$ may also be 5.

Rev. 18'-20': There are traces of signs in 18', but 19'-20' may have been uninscribed.

This large and very thick fragment (9.5 x 14.0 x 2.6-5.0 cm) preserves the upper (lower) half of the obverse (reverse) of a tablet. No physical join with other fragments could be established. A portion of unknown height is missing. The lower half of the obverse and most of the reverse are heavily eroded and almost illegible. The obverse contains a table with regular numbers and their reciprocals which might continue at the top of the reverse. Text C is the only known tablet from Babylon preserving Nrs. 1-3 of the standard table and the second fragment, along with Text D, to preserve Nr. 4. No satisfying interpretation of the reverse could be found. There appears to be a single column, but even that is not certain. For most of the badly eroded signs no plausible reading could be established. Near the bottom there are sections delimited by horizontal rulings. In Rev. 16' the regular number 2.1.4.8.3.0.27, equivalent to $3^{23}$, is written at some distance from the left edge. This number is mentioned or implied in several other tablets,[21] but its significance in Text C is unclear. None of the other numbers on the reverse could be identified.

## Text D = BM 37095 (80-6-17, 844)

Obverse

| | i = *s* | ii = 1/*s* | Nr |
|---|---|---|---|
| 1' | [IGI 1.1.30.33.4]⸢5⸣ | ⸢58.3⸣[1.39.35.18.31.6.40] | 4 |
| 2' | [IGI 1.1.43.42.13].⸢20⸣ | 5⸢8.2⸣9.[12] (error) | 5 |
| 3' | [IGI 1.2.12.28.4]⸢8⸣ | 57.52.1[3.20] | 6 |
| 4' | [IGI 1.2].⸢30⸣ | 57.⸢3⸣[6] | 7 |
| 5' | [IGI 1.3.12.35.33.20] | 5⸢6⸣.5⸢7.11⸣.[15] | 8 |

Critical commentary:
Obv. 2': 29 is an error for 19.

This fragment (6.0 x 2.5 x 2.6 cm) is inscribed on one side which is most likely the obverse; the other side is destroyed. No edges of the tablet are preserved. Since the numbers belong to the beginning of the standard table, not much clay may be missing above the fragment. It might belong to the same tablet as Text E, but there is no physical join.

## Text E = BM 33447 (Rm4, 1)

Obverse

| | i = *s* | ii = 1/*s* | Nr |
|---|---|---|---|
| 1' | [IGI 1.5.55.4.41].15 | [54.36.48] | 15 |
| 2' | [IGI 1.6.21.18].⸢4⸣3.12 : 1.6.⸢40⸣ | [54.15.12.30 : 54] | 16, 17 |
| 3' | [IGI 1.7.30] | 5⸢3.20⸣ | 18 |
| 4' | [IGI 1.8.1]⸢6⸣ | 52.4⸢4⸣.[3.45] | 19 |
| 5' | [IGI 1.8].⸢20⸣.37.30 | 52.⸢40⸣.[29.37.46.40] | 20 |
| 6' | [IGI 1.9].⸢7⸣.12 | 52.⸢5⸣ | 21 |
| 7' | [IGI 1.9.2]⸢6⸣.40 | 51.⸢50⸣.[24] | 22 |
| 8' | [IGI 1.9.5]⸢9⸣.2.24 | 51.26.[25.11.6.40] | 23 |
| 9' | [IGI 1.10.1]⸢8⸣.45 | 5⸢1⸣.1[2] | 24 |
| 10' | [IGI 1.11].6.40 | 50.⸢3⸣[7].⸢30⸣ | 25 |
| 11' | [IGI 1.11.11].29.3.45 | 50.34.4.2[6.40] | 26 |
| 12' | [IGI 1.12] | 50 | 27 |
| 13' | [IGI 1.12.49].⸢4⸣ | 49.1⸢5?⸣.18.19.[...] (error) | 28? |
| 14' | [IGI 1.12.54] | 49.22.57.⸢4⸣[6.40] | 29 |
| 15' | [IGI 1.13.9.34.29.8.8.5]⸢3.20 | 49.12.27⸣ | 30 |

21 It is present in the extended table of reciprocals *TU* 31, from Uruk. Its fourth power, $9^{46}$, appears in Text A, Obv. 2, and its square has been tentatively identified in Text K, X.5'.

Reverse

| | | | |
|---|---|---|---|
| 1' | [IGI 1.38.18.14.24] | 3⸢6.37⸣.[15.56.15] | 70 |
| 2' | [IGI 1.38.24.54] | 36.3⸢4⸣.[47.14.34.4.26.40] | 71 |
| 3' | [IGI 1.38.45.55.3]3.20 | 36.[27] | 72 |
| 4' | [IGI 1.38.52.37.1].⸢5⸣2.30 | 3⸢6⸣.[24.32] | 73 |
| 5' | [IGI 1.40] | ⸢3⸣[6] | 74 |
| 6' | [IGI 1.41.8.8.5]⸢3⸣.20 | [35.35.44.31.52.30] | 75 |
| 7' | [IGI 1.41.15 | 35.33.20] | 76? |
| 8' | [IGI 1.42.24 | 35.9.22.30] | 77? |
| 9' | [IGI 1.42.52.50.2]⸢2⸣.13.[20 | 34.59.31.12] | 78 |

Critical commentary:
Obv. 2': This line contains two pairs of numbers.
Obv. 13': The 5 might also be an 8. One expects 49.26.18.30.56.15. The origin of this error is not clear.

This fragment (5.9 x 6.6 x 3.5 cm) does not preserve any edge of the tablet. Each side probably contained one half of the standard table. The reverse might physically join the fragment Liverpool 29.11.77.34 (Friberg 2005).[22]

## Text F = BM 32681 (76-11-17, 2450)

Obverse

| | i = *s* | ii = 1/*s* | Nr |
|---|---|---|---|
| 1' | [IGI 1.11.6.40 | 50.3]⸢7.30⸣ | 25 |
| 2' | [IGI 1.11.11.29.3.45 | 50.3]⸢4⸣.4.2⸢6.40⸣ | 26 |
| 3' | [IGI 1.12 | 50] | 27 |
| 4' | [IGI 1.12.49.4 | 49.2]⸢6.18.30.5⸣6.15 | 28 |
| 5' | [IGI 1.12.54 | 49].⸢2⸣2.⸢5⸣[7.4]⸢6.40⸣ | 29 |
| 6' | [IGI 1.13.9.34.29.8.8.53.20 | 49].⸢1⸣2.⸢2⸣[7] | 30 |
| 7' | [IGI 1.13.14.31.52.30 | 49.9.7].1⸢2⸣ | 31 |
| 8' | [IGI 1.13.43.40.48 | 48.49.41].⸢15⸣ | 32 |
| 9' | [IGI 1.14.4.26.40 | 48.36] | 33 |
| 10' | [IGI 1.14.38.58.33.36 | 48.13].⸢31.6.40⸣ | 34 |
| Reverse | | | |
| 1' | [IGI 1.34.48.53.20 | 37.58].⸢7⸣.30 | 64 |
| 2' | [IGI 1.34.55.18.45 | 37.55].33.20 | 65 |
| 3' | [IGI 1.36 | 37].⸢30⸣ | 66 |
| 4' | [IGI 1.36.27.2.13.20 | 37.1]⸢9⸣.29.16.48 | 67 |
| 5' | [IGI 1.37.12 | 37.2.1]⸢3⸣.20 | 68 |
| 6' | [IGI 1.37.39.22.30 | 36.51].⸢50⸣.24 | 69 |
| 7' | [IGI 1.38.18.14.24 | 36.37.15].⸢56⸣.15 | 70 |
| 8' | [IGI 1.38.24.54 | 36.34.47.1]⸢4⸣.34.4.26.⸢40⸣ | 71 |
| 9' | [IGI 1.38.45.55.33.20 | 36.2]⸢7⸣ | 72 |
| 10' | [IGI 1.38.52.37.1.52.30 | 36.24].32 | 73 |

This small fragment (2.2 x 5.0 x 1.6-2.3 cm) from the right edge of a tablet preserves several final digits of reciprocal numbers. Each side of the tablet must have contained one half of the standard table. No physical join with other

22 Since the Liverpool fragment is lost, this cannot be verified. No physical join exists with candidates BM 34577 (Vaiman 1961, Text F), BM 76984 (Britton 1991-3, Text A) or BM 77051 (Britton 1991-3, Text B).

fragments could be established.[23]

## Text G = BM 42980 (81-7-4, 744)

Reverse

| | i = *s* | ii = 1/*s* | Nr |
|---|---|---|---|
| 1 | [IGI 1.23.20 | 4]⸢3⸣.12 [...] | 49 |
| 2 | [IGI 1.24.22.30 | 42].⸢40⸣ [...] | 50 |
| 3 | [IGI 1.25.20 | 42].⸢1⸣1.15 [...] | 51 |
| 4 | [IGI 1.26.24 | 41.40] [...] | 52 |
| 5 | [IGI 1.26.48.20 | 41.2]⸢8⸣.19.12 [...] | 53 |
| 6 | [IGI 1.27.47.29.22.57.46.40 | 41.0].22.30 [...] | 55 |
| 7 | [IGI 1.27.53.26.15 | 40.57].⸢3⸣6 [...] | 56 |
| 8 | [IGI 1.28.53.20 | 40.30] [...] | 57 |
| 9 | [IGI 1.31.7.30 | 39.30.22].⸢1⸣3.⸢20⸣ [...] | 59 |
| 10 | [IGI 1.32.9.36 | 39.3.4]⸢5⸣ [...] | 60 |

This fragment (3.5 x 4.0 x 1.6 cm) is inscribed on one side which is most likely the reverse; the obverse is destroyed. A small segment of the upper edge is preserved. The numbers match the beginning of the second half of the standard table. No physical join with other fragments could be established. Two entries, Nrs. 54 and 58, appear to be lacking, but perhaps some lines contained two pairs, as in Text E.

## Text H = BM 36917 (80-6-17, 658)

Side X

| | i' = *s* | ii' = 1/*s* | iii' = *s* | iv' = 1/*s* |
|---|---|---|---|---|
| 1' | [xxxx | x | IGI.25].⸢GAL$_2$.BI⸣ | [2.24] |
| 2' | [xxxx | x] | IGI.27.⸢GAL$_2$.BI | 2⸣.[13.20] |
| 3' | [xxxx | x] | IGI.30.GAL$_2$.BI | [2] |
| 4' | [IGI.4.GAL$_2$.BI | 1]⸢5⸣ | IGI.32.GAL$_2$.BI | ⸢1⸣.[52.30] |
| 5' | [IGI.5.GAL$_2$.BI | 1]⸢2⸣ | IGI.36.GAL$_2$.BI | [1.40] |
| 6' | [IGI.6.GAL$_2$.BI | 10] | IGI.40.GAL$_2$.BI | [1.30] |
| 7' | [IGI.8.GAL$_2$.BI] | ⸢7⸣.30 | IGI.45.⸢GAL$_2$⸣.[BI | 1.20] |
| 8' | [IGI.9.GAL$_2$.BI] | ⸢6.40⸣ | IGI.48.⸢GAL$_2$⸣.[BI | 1.15] |
| 9' | [IGI.10.GAL$_2$.BI | 6] | ⸢IGI.50⸣.[GAL$_2$.BI | 1.12] |

This fragment (5.5 x 4.5 x 2.5 cm) is inscribed on one side; the other side is destroyed. It does not preserve any edge of the tablet. The preserved text matches the standard OB table of 30 regular numbers and their reciprocals (Neugebauer & Sachs 1945: p. 11). This suggests that the table begins in X.1' and that there were seven more lines below X.9'. The tablet may be a duplicate of BM 34592 (Aaboe 1965, Text IV), a LB fragment from Babylon on which the same table is followed by an OB-style table of squares.

## III. Tables of squares of regular numbers

Three newly identified fragments contain squares of regular numbers. All are copies or variants of a table with squares of roughly the same one hundred regular numbers known from the tables of reciprocals (Section II). This brings the total number of such fragments to ten.[24] Text K includes squares which are not in the standard table. A peculiar feature shared by all known LB fragments of tables of squares is that none preserves a column with the principle numbers (*s*). It

23 A partly overlapping stretch is preserved in Text E, but they do not belong to the same tablet.

24 Previously published texts: Vaiman (1961), Texts G-J (=Aaboe 1965, Texts V-VIII); Aaboe (1965), Texts II-III (Text J published here joins the latter); *LBAT* 1640 (Friberg 1986, p. 87); Britton (1991-3), Text D.

is therefore likely that these tables lack such a column,[25] which raises the interesting question of how they were used in practice.[26]

## Text I = BM 37020 (80-6-17, 764)

Obverse?

| | $s^2$ | $s$ | Nr |
|---|---|---|---|
| 1' | [xxxxxxxxxxxxxxxxxx]⸢x⸣[...] | 1.... | ? |
| 2' | [xxxxxxxxxxx]⸢x⸣[xxx].4⸢5⸣.[...?] | 1.... | ? |
| 3' | [1.28.22.25.4]⸢3⸣.3⸢2⸣.1⸢6⸣ [:] 1.28.34.2⸢4⸣.[36] | 1.12.49.4 1.12.54 | 28, 29 |
| 4' | [1.29.12.19.2]⸢6⸣.34.23.19.49.⸢3⸣[8.8.36.52.20.44.26.40] | 1.13.9.34.29.8.8.53.20 | 30 |
| 5' | [1.29.24.25].4.54.26.0.56.1⸢5⸣ [...] | 1.13.14.31.52.30 | 31 |
| 6' | [1.30.35.49].4.44.32.38.24 [...] | 1.13.43.40.48 | 32 |
| 7' | [1.31.26.58].⸢6⸣.2⸢5⸣.11.6.40 : ⸢1.30⸣[...] | 1.14.4.26.40 1.... | 33, ? |
| 8' | [xxxxxx : 1.3]⸢5.5⸣3.30.12.20.⸢4⸣[4.26.40 ...] | 1.... 1.15.51.6.40 | ?, 36 |
| 9' | [1.36.6.30.14.3].4⸢5⸣ : ⸢1⸣.[...] | 1.15.56.15 1.... | 37, ? |
| 10' | [xxxxxxxx]⸢xx⸣[...] | 1.... | ? |

Critical commentary:
Obv. 1': There are traces of a digit 10-50.
Obv. 2': The preserved digit might belong to 1.24.28.12.58.45.57.7.44.3.45 (Nr. 26).
Obv. 7': The second square may be 1.32.52.33.46.0.30.52.24.57.36 (Nr. 34) or 1.33.45 (Nr. 35).
Obv. 8': The available space at the beginning suggests that this line contained another square, perhaps 1.33.45 (Nr. 35).
Obv. 9': The second square may be 1.38.18.14.24 (Nr. 38).
Obv. 10': There are traces of a 10 (or 20-50) followed by a 4-8.

This fragment (4.9 x 4.7 x 1.2 cm) is inscribed on one side which is probably the obverse; the other side is destroyed. It does not preserve any edges of the tablet. No physical join with other tables of squares could be established. Some lines contain two squares, a feature known from other tablets.[27]

## Text J = BM 32178 (76-11-17, 1905+2228) +Rm 848

Obverse

| | $s^2$ | $s$ | Nr |
|---|---|---|---|
| 1' | [1.26.2]⸢4⸣ | 1.12 | 27 |
| 2' | [1.28.22].25.⸢43.32⸣.[16] | 1.12.49.4 | 28 |
| 3' | [1.28.3]⸢4.2⸣4.36 | 1.12.54 | 29 |
| 4' | [1.29.12.19].26.34.23.19.49.[38.8.36.52.20.44.26.40] | 1.13.9.34.29.8.8.53.20 | 30 |
| 5' | [1.29.24.25].4.54.26.0.5⸢6⸣.[15] | 1.13.14.31.52.30 | 31 |
| 6' | [1.30.35.49].⸢4⸣.44.32.38.[24] | 1.13.43.40.48 | 32 |
| 7' | [1.31.26.58.6.25].⸢11.6.40⸣ | 1.14.4.26.40 | 33 |

The fragment Rm 848 (5.9 x 4.9 x 1.5 cm) physically joins BM 32178 (Aaboe 1965, Text III) at the top of the obverse. Its reverse is destroyed. The joined fragments (9.8 x 14.5 x 2.6-3.6 cm) preserve Nrs. 27-54 (obverse) and Nrs. 55-79 (reverse) of the standard table. Hence the tablet must have contained the entire table. Only Rm 848 is edited here; for

25 This is certain for at least one fragment listing squares in two adjacent columns (Vaiman 1961, Text G, probably also Text H).

26 Cf. also Vaiman (1961), p. 64; Aaboe (1965), p. 86. In principle a square could be looked up by placing the table of squares next to a standard table of reciprocals, but occasional differences in layout cast some doubt on the practicality of this procedure.

27 Aaboe (1965), Texts III, VIII; Britton (1991-3), Text D. In other tables, e.g. Aaboe (1965), Text VII, this feature is clearly lacking.

BM 32178 cf. Aaboe (1965).

## Text K = BM 45884 (81-7-6, 315)

Side X

| | $s^2$ | $s$ | Nr |
|---|---|---|---|
| 1' | [3.48.52].⸢54.36.33⸣.4⸢5⸣ | 1.57.11.15 | 98 |
| 2' | [3.51.55].⸢41.38⸣.32.25.[57.30.14.24] | 1.57.57.53.16.48 | 99 |
| 3' | [3.52.27.8.15.22.27].⸢50⸣.24 | 1.58.5.52.48 | - |
| 4' | [3.54.6.38.21].14.4.26.⸢40⸣ | 1.58.31.6.40 | 100 |
| 5' | [4.4.17.40.45].⸢1⸣3.1⸢7.45⸣.[52.14.42.12.9] | 2.1.4.8.3.0.27 | - |

Critical commentary:
X.1': ⸢54.36.33⸣.4⸢5⸣ replaces 21.33.32.4[4...] (*ACT*).
X.2': ⸢41⸣ replaces 21 (*ACT*).
X.3': ⸢50⸣ replaces ...1]4 (*ACT*).
X.5': The reconstruction of this number is tentative but plausible. 1⸢7⸣ replaces 18 (*ACT*).

This small fragment (4.7 x 4.5 x 2.5 cm) does not preserve any edge of the tablet. It is inscribed on one side; the other side is destroyed. It was previously published, in transliteration only, as *ACT* 1001.[28] No physical join with other tables of squares could be established. The fragment is probably located near the lower edge of the reverse. The only other fragment preserving Nrs. 98, 99 and 100 is *LBAT* 1640 (Friberg 1986, p. 85). The squares reconstructed in X.3' and 5' are attested here for the first time; they are not in the standard table.

# IV. Other tables

The two remaining fragments do not form a homogeneous group. Both contain multiplications of one kind or another.

## Text L (BM unnumbered fragment 3.6): combined multiplication table

Side X

| | i' = *n* | ii' = *n*·52.30 | iii' = *n* | iv' = *n*·? |
|---|---|---|---|---|
| 1 | [... 1] | 5⸢2.30⸣ | 1 | ⸢30⸣[...] |
| 2 | [... 2] | 1.45 | 2 | 1.[...] |
| 3 | [... 3] | ⸢2⸣.37.30 | 3 | 1.⸢30$^{?}$⸣[...] |
| 4 | [... 4] | ⸢3.30⸣ | 4 | 2.1⸢2$^{?}$⸣[...] |
| 5 | [... 5 | 4.22].⸢30⸣ | 5 | 2.⸢40$^{?}$⸣[...] |
| 6 | [... 6 | 5.15] | ⸢6 | 3⸣.[...] |

Critical commentary:
X.iv'.3: The digit following 1 begins with 30-50.
X.iv'.4: The numeral following 2.10 is 2-8.
X.iv'.5: The digit following 2 begins with 30-50.

This small fragment (4.5 x 2.0 cm) from the upper edge of a tablet has not received a registration number.[29] It is inscribed on one side; the other side is destroyed. It is probably a combined multiplication table. Very common in the OB era, these tables list products of various principle numbers with single-digit factors increasing from 1 to 20 in steps of 1 and from 30 to 50 in steps of 10 (Neugebauer & Sachs 1945: pp. 24-33). The principle number in Col. ii' is 52.30; the digits in iv'.4 imply that the principle number in this column is between 33 (= 2.12 / 4) and 34.30 (= 2.18 / 4). Only two LB fragments of a combined multiplication table were known up to now: U 91 from Uruk (Aaboe 1969) and BM 36849 from Babylon (Aaboe 1999). Both are partly overlapping duplicates of presumably the same table. The present

28 Neugebauer (*ACT*) assumed an astronomical content but he did not identify the numbers.
29 Numerous small, unnumbered fragments from the 'Babylon' collection thought to have an astronomical content are kept in five plastic boxes (as of 2004); this fragment is kept in box 3.

fragment is different from the known OB and LB types, so it must belong to a hitherto unknown variant.

## Text M = BM 37338 (80-6-17, 1095)

| | 1.32.3.21.47.51.6.40·*f* | | *f* |
|---|---|---|---|
| Side X | | | |
| 1' | 1.⸢4⸣.26.21.⸢15⸣.[29.46.40 | ...] | 42 |
| 2' | 1.5.58.24.3⸢7.17⸣.[37.46.40 | ...] | 43 |
| 3' | 1.7.30.27.59.5.2⸢8⸣.[53.20 | ...] | 44 |
| 4' | 1.9.2.31.20.5[3.20 | ...] | 45 |
| 5' | ⸢1.10.34.34.42⸣.[41.11.6.40 | ...] | 46 |
| Side Y | | | |
| 1' | ⸢11.39.37.3⸣3.39.[40.26.40 | ...] | 7.36 |
| 2' | 11.58.2.14.1.14.40 [ | ...] | 7.48 |
| 3' | 12.16.26.54.22.4[8.53.20 | ...] | 8 |
| 4' | 12.3⸢4.51.34⸣.[44.23.6.40 | ...] | 8.12 |

This fragment (5.5 x 4.2 x 2.0-2.5 cm) includes a segment of the left edge of a tablet; no other edges are preserved. On both sides one column of numbers is partly preserved. Only the numbers in X.4' and Y.3' were identified in a computer-generated table of all regular numbers with up to thirty digits, implying that the other numbers are not regular. Hence Text M is not a table of reciprocals or powers of regular numbers. By querying a computer-generated table of products of all integers 1...59 with all regular numbers with up to twenty digits, the preserved digits matched products of 1.32.3.21.47.51.6.40 (=$2^{35}\cdot5^{3}$) with two sequences of linearly increasing small factors (*f*). On side X the line-by-line increment of *f* is 1; on side Y it is 12, resulting in a line-by-line difference 12 · 1.32.3.21.47.51.6.40 = 18.24.40.21.34.13.20. If side X covered *f*=1...59 then forty-one additional lines must be reconstructed above X.1' and 13 below X.5', but that does not yield a plausible reconstruction, so the full extent of Text M remains unknown. The significance of the numbers is also unclear. Linearly changing functions are common in Babylonian mathematical astronomy, but a plausible astronomical interpretation could not be found.

## Acknowledgements

I thank the Trustees of the British Museum for providing access to their cuneiform collection and for permission to publish the texts.

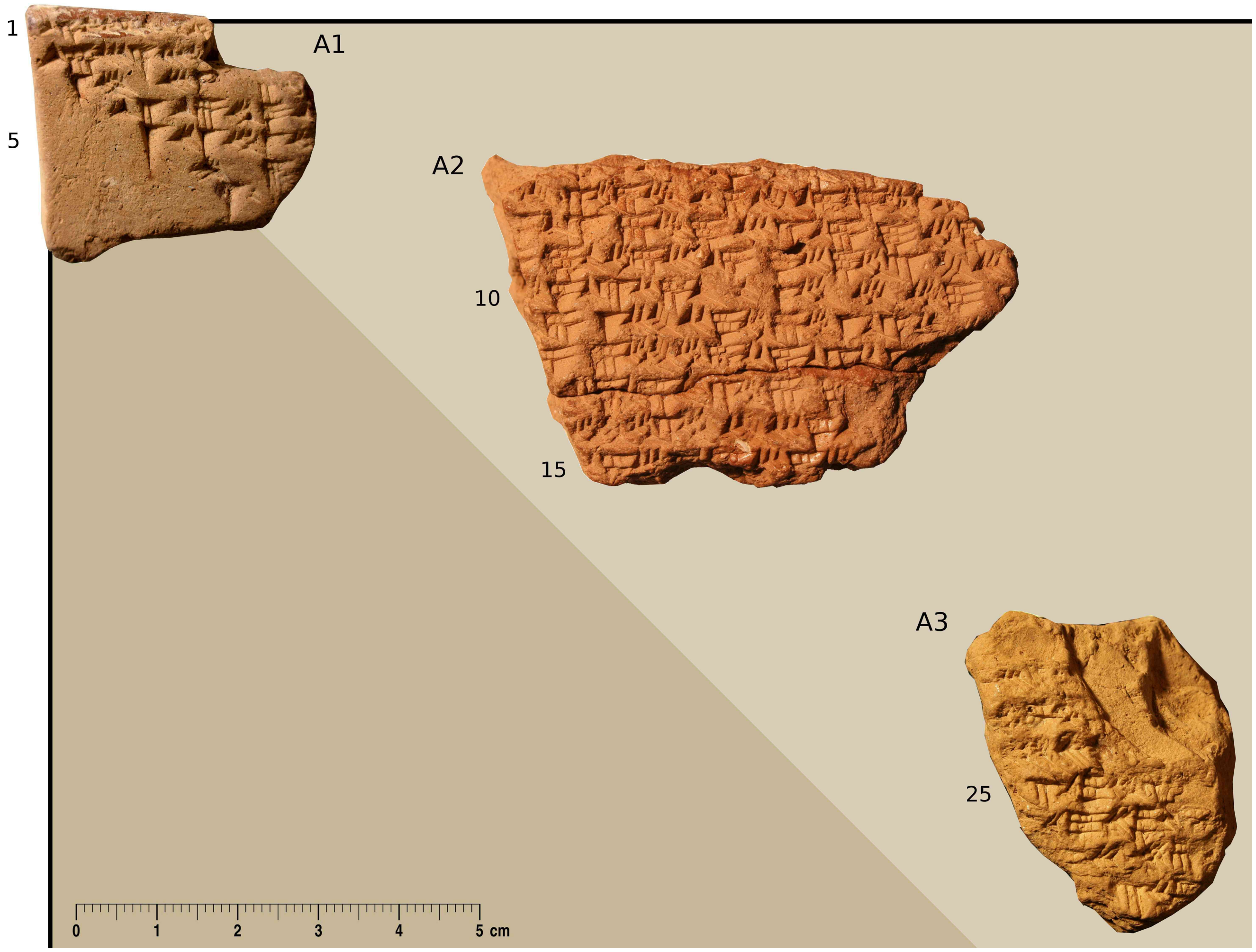

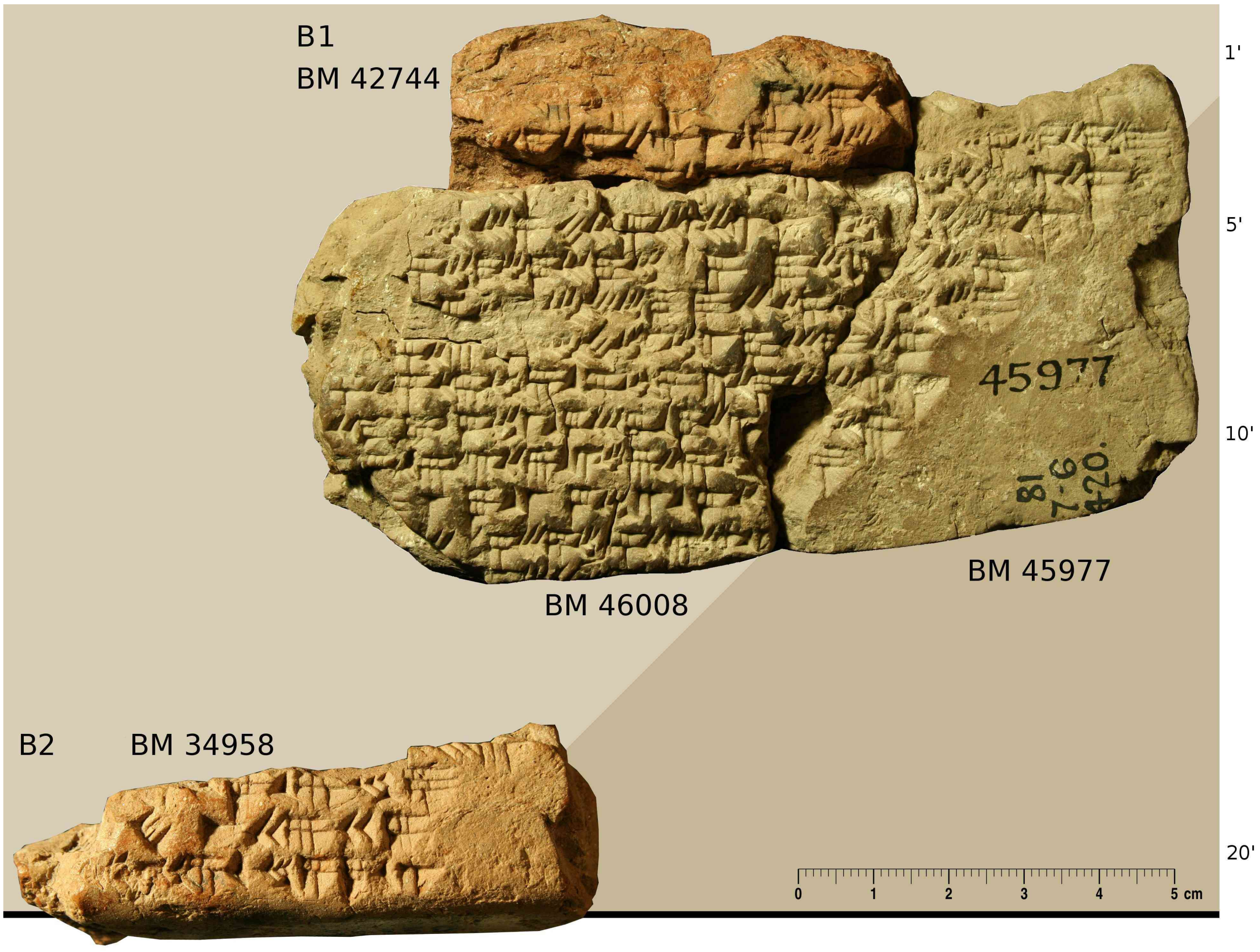

C

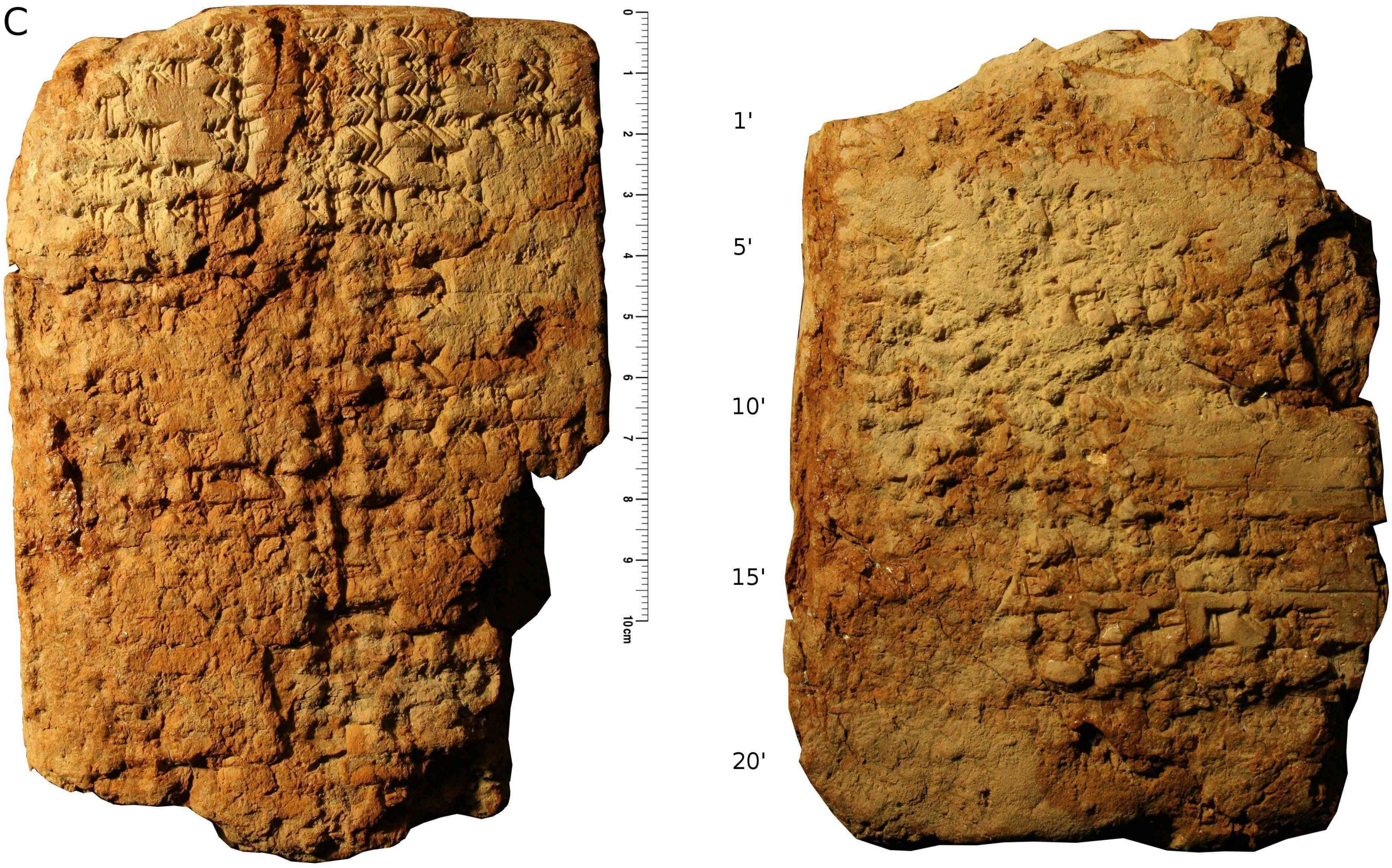

D

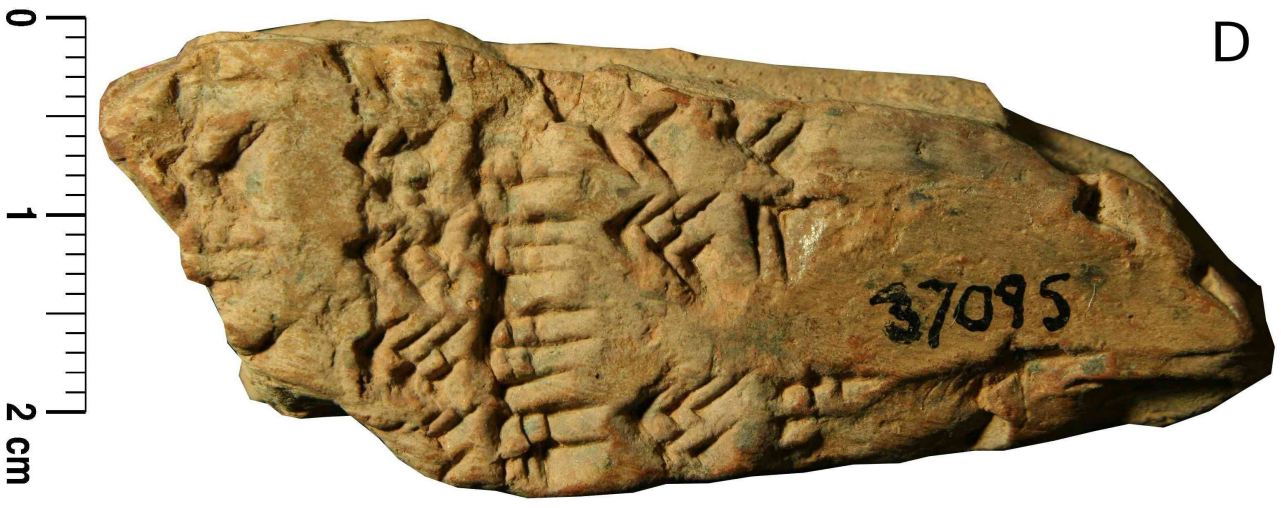

L

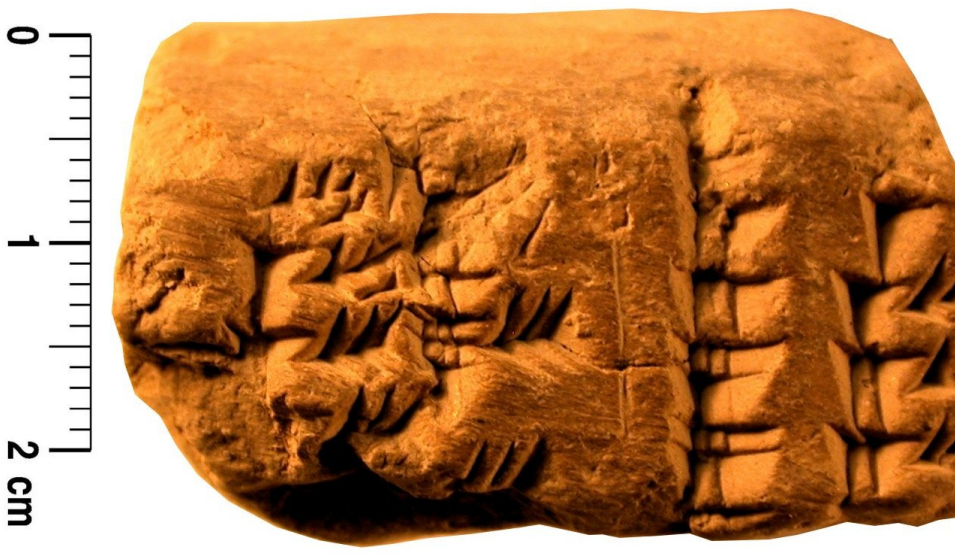

E

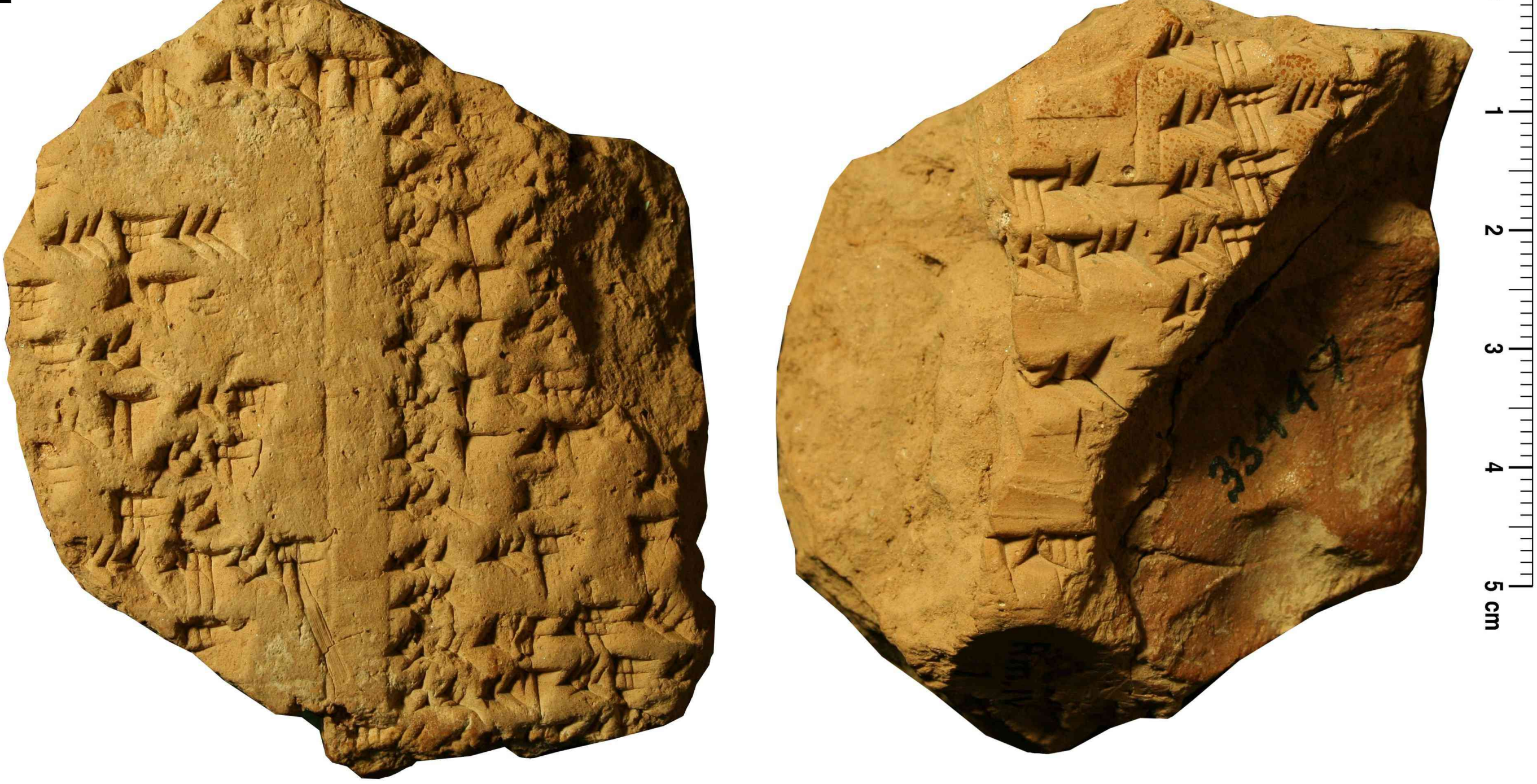

F

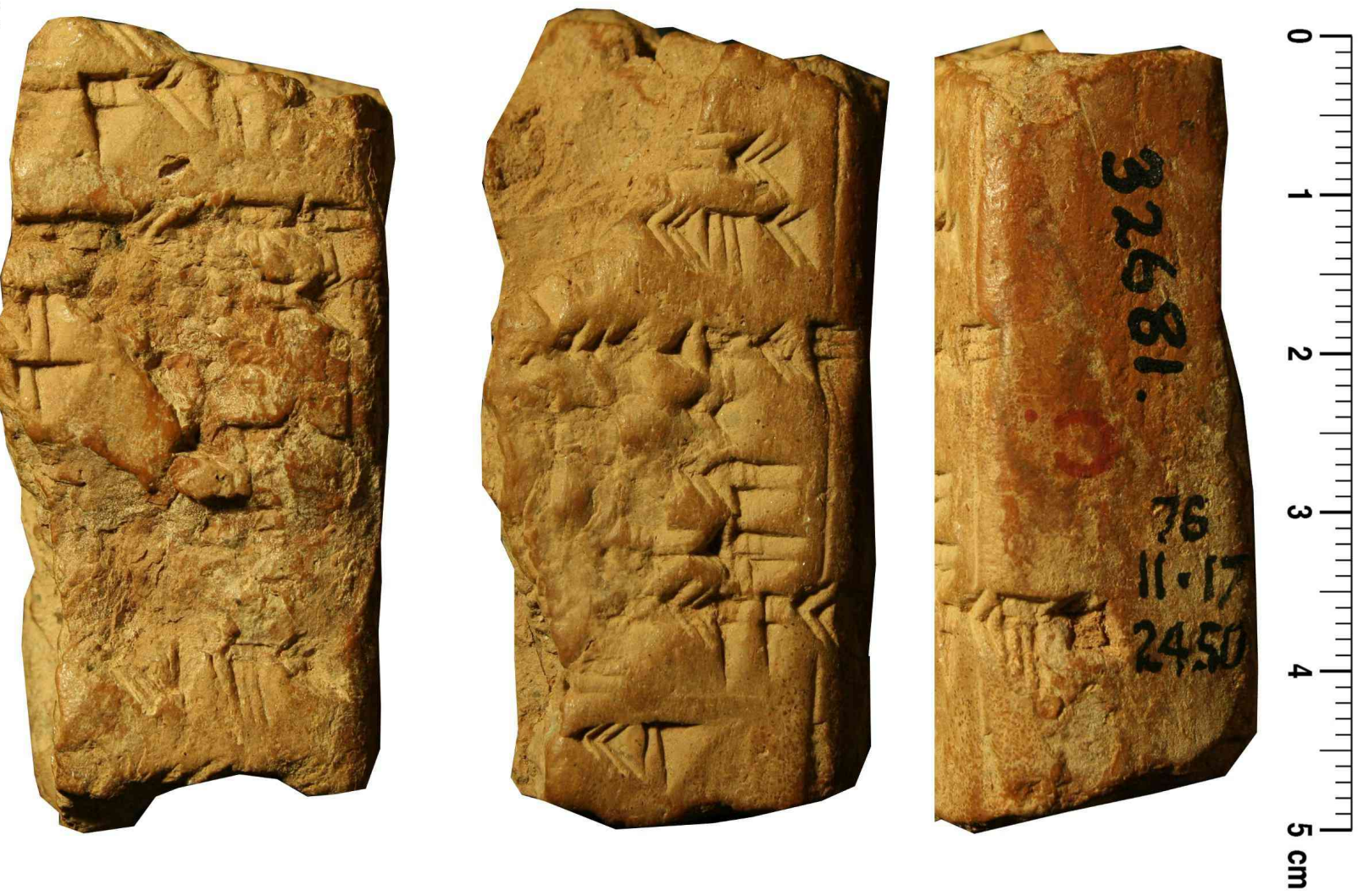

G

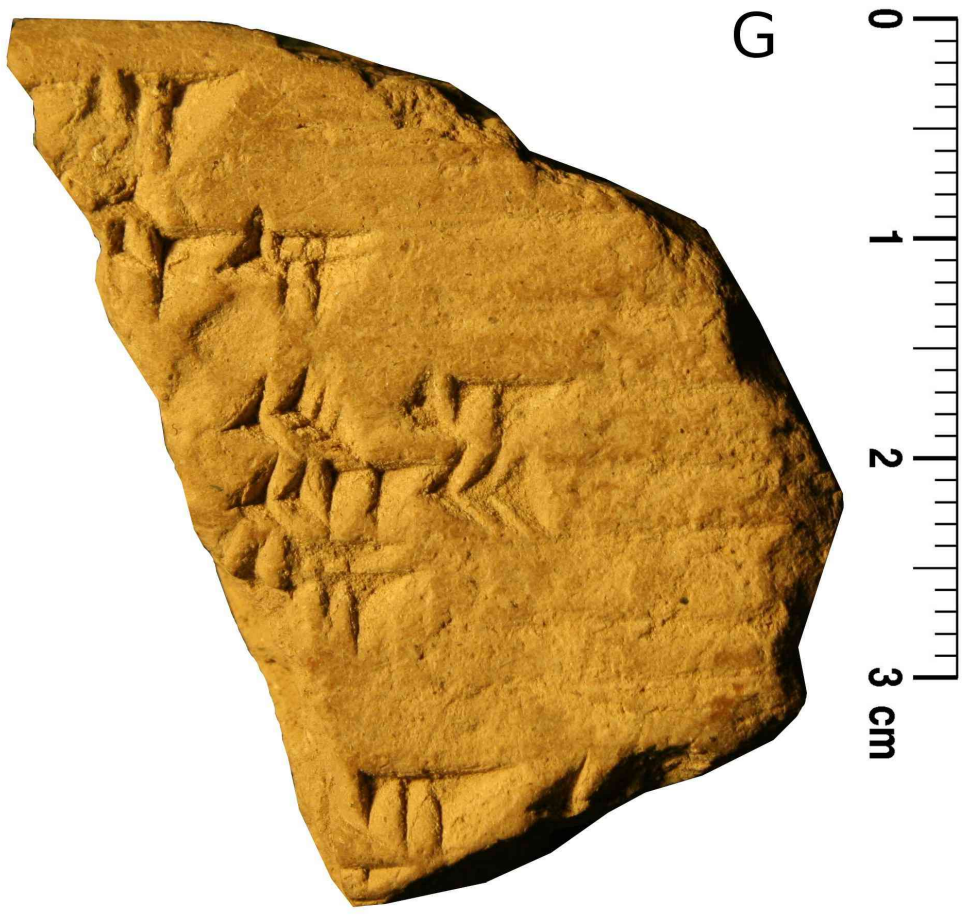

H

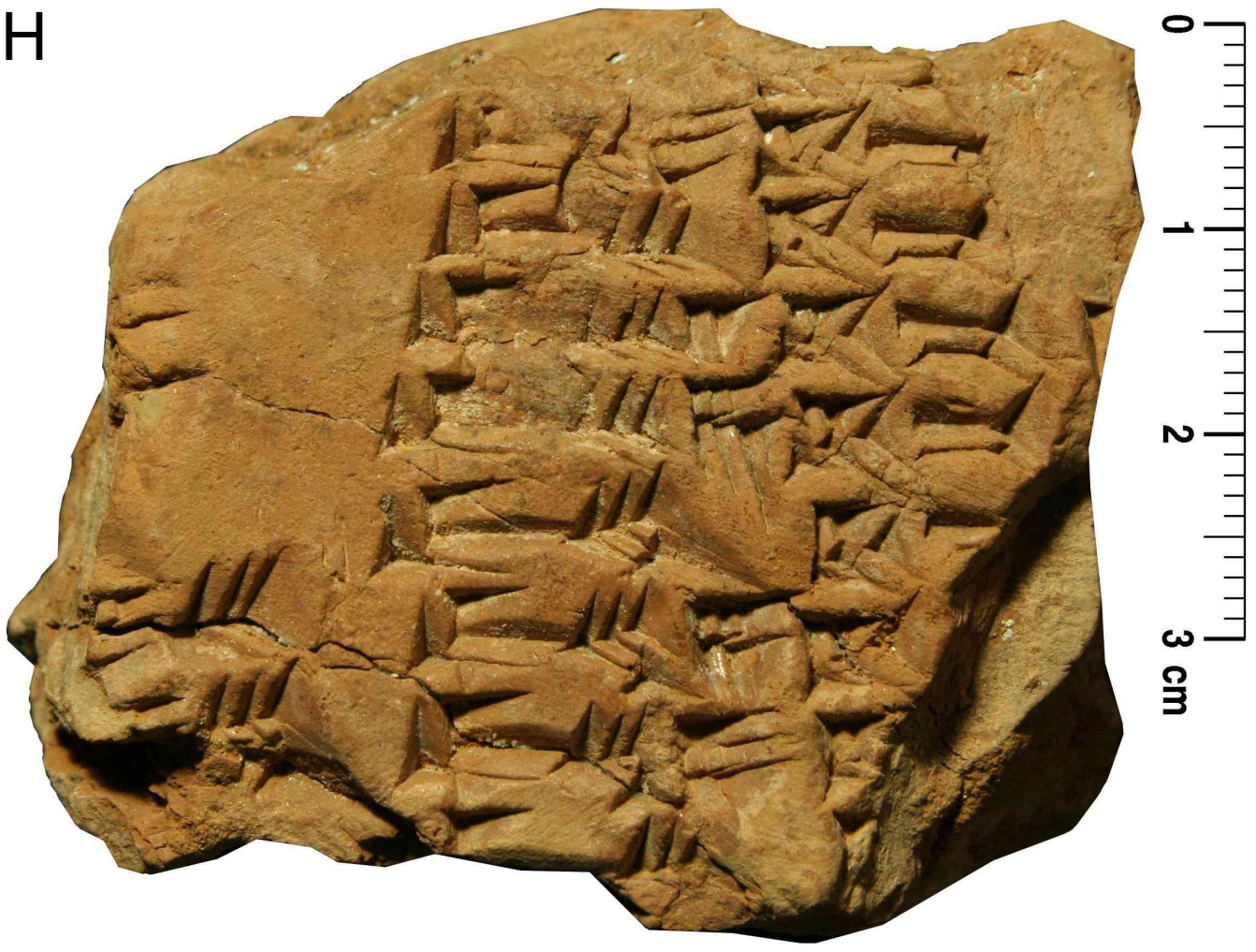

I

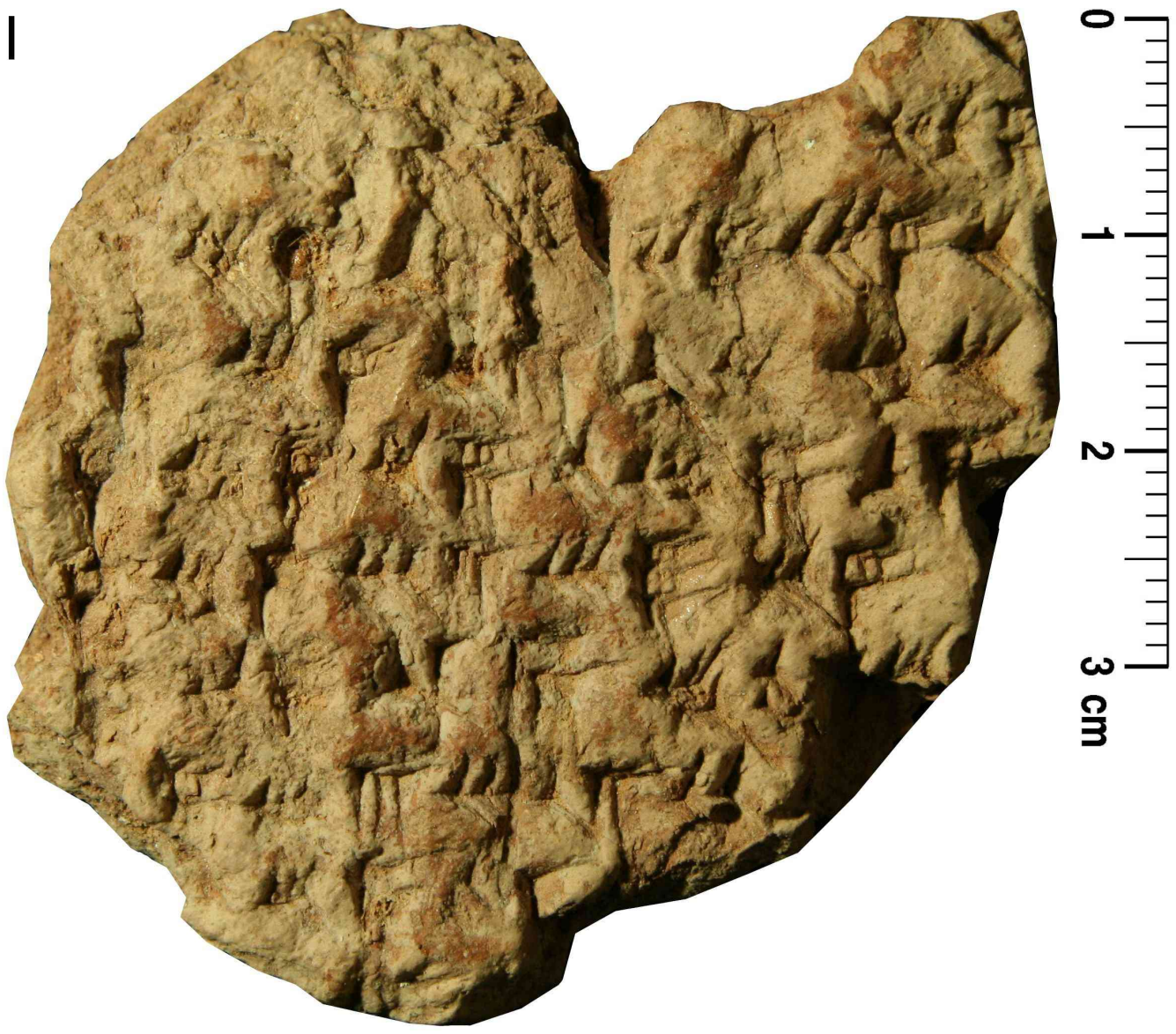

J

Rm 848

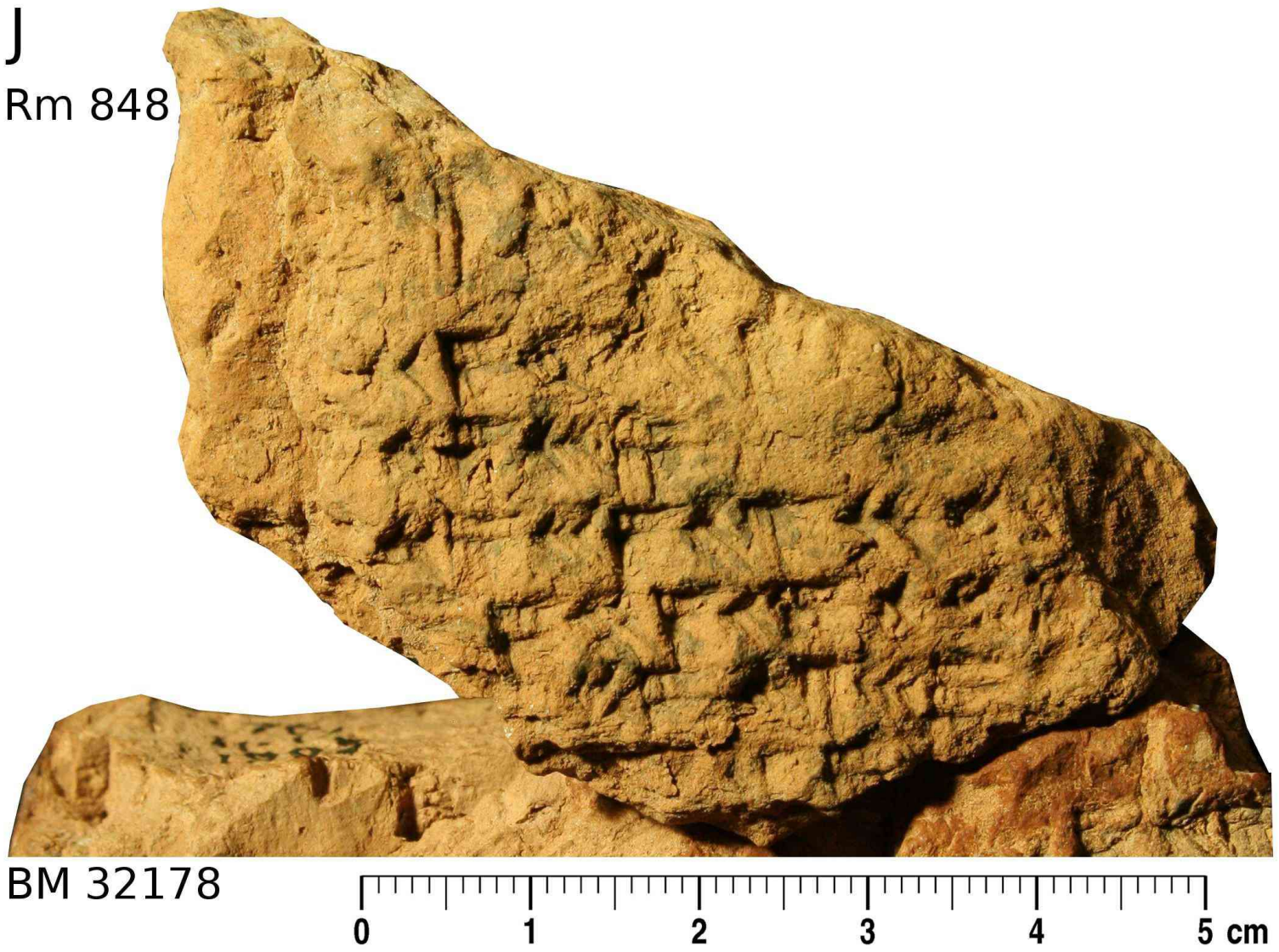

BM 32178

K

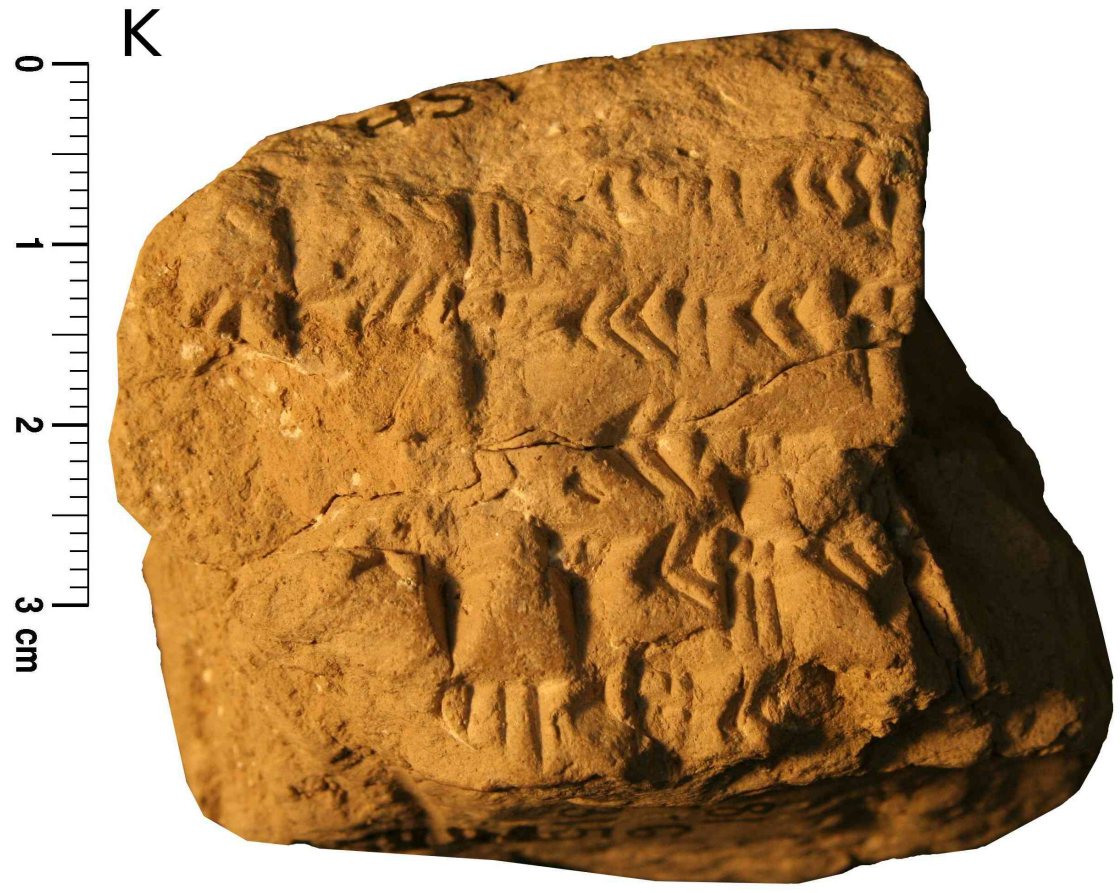

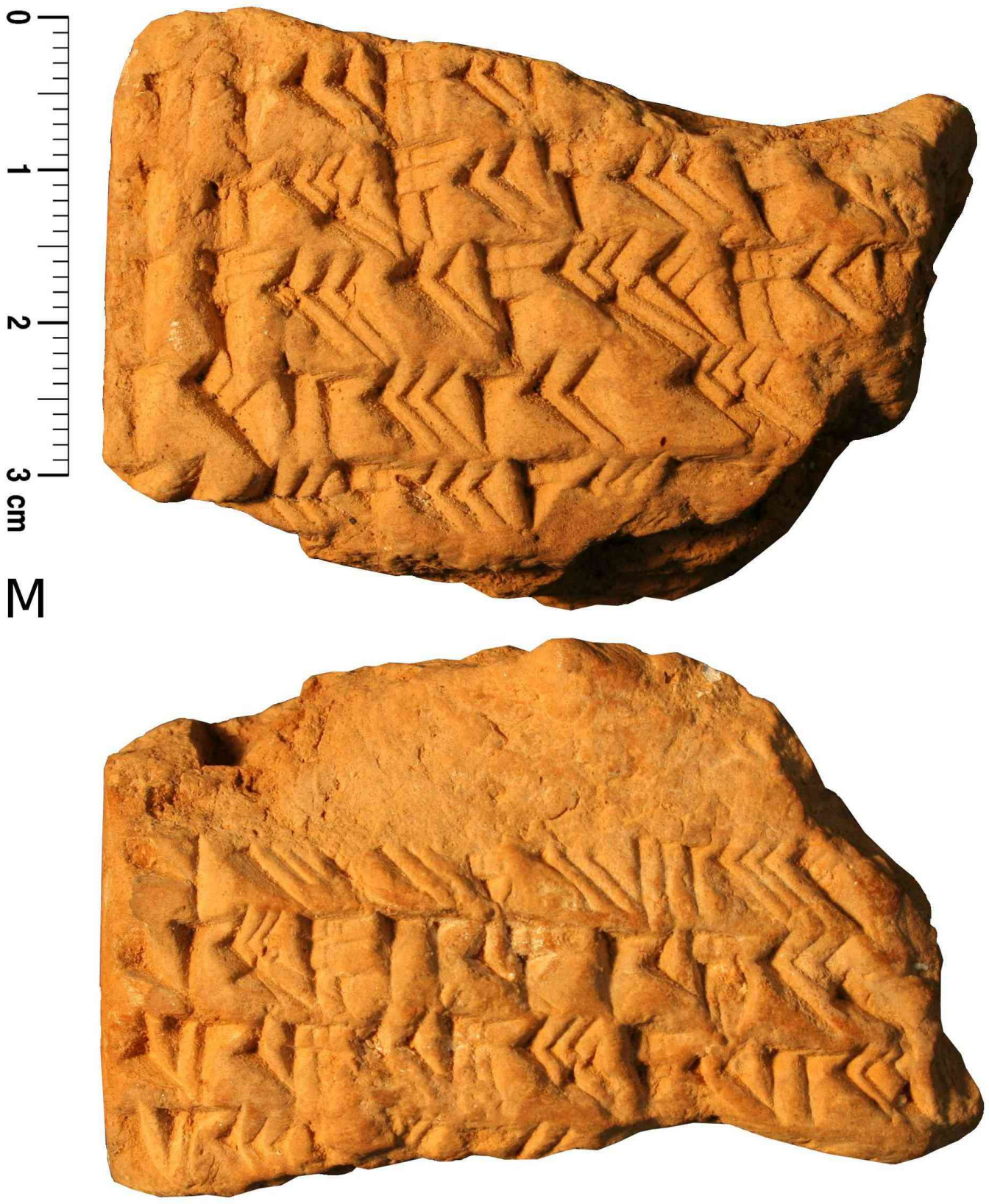

M